\documentclass[12pt]{article}
\usepackage[margin = 1in]{geometry}
\usepackage{amsfonts, amsmath, amssymb}
\usepackage[none]{hyphenat}
\usepackage{fancyhdr}
\usepackage{amsmath}
\usepackage{xcolor}
\usepackage{graphicx}
\usepackage{amsthm}
\usepackage{amssymb}
\usepackage{caption}
\usepackage{subcaption}
\usepackage{float}
\usepackage{enumitem}
\usepackage{blindtext}
\usepackage{lipsum}
\usepackage{hyperref}

\title{How many ways to color the map of America?}
\date{\vspace{-5ex}}

\author{Rebekah Bessett$^{\dagger}$, Jennifer Canizales${^\dagger}$, Jasbir S. Chahal$^{*}$,\\
Thomas Fackrell$^{\dagger}$, Vanessa Rico$^{\dagger}$}

\begin{document}
\maketitle
\renewcommand{\headrulewidth}{0pt} 
\medskip
\begin{abstract}
    Although the Four Color Conjecture originated in cartography, surprisingly, there is nothing in the literature on the number of ways to color an actual geographic map with four or less colors. In this paper, we compute these numbers, with exponentially increasing order of difficulty, for the maps of Canada, France, and the USA. Our attempts to compute the latter two lead to some new results on the chromatic polynomial of graphs.
    \let\thefootnote\relax\footnotetext{
    $* $ Faculty Advisor\\
    $\dagger$ Undergraduate Students
    } 
\end{abstract}

\section{Introduction}
The \textit{Four Color Theorem} conjectured by Francis Guthrie in 1852, asserts that no planar graph needs more than four colors to color it properly. A coloring is proper if no two vertices which are adjacent are assigned the same color. A related problem, on which we found no literature for real life graphs, is to determine the number of ways to color the graph $G_X$ of a given country X, with $k(G_X)\leq 4$ colors, $k(G_X)$ being the smallest number of colors with which $G_X$ can be colored properly. It is easy to see that $k(G_X) = 3$ if $G_X$ is the graph of the map of the provinces/territories of Canada. However, a Google search  found no entry for the number of ways to color the map of Canada with three colors. As we will show this number is 576 and fairly easy to compute by hand. But this is not so easy for the map of the contiguous regions of France. There are programs which use the computing power of a machine that will find this number for the main part of France which has twelve contiguous regions. In fact, we used different programs to compute and compare this number for the map of the 29 counties of the state Utah of USA. However, these programs reach their limit when the number of the vertices of a graph is around 40. Although the Four Color Conjecture originated in Cartography, surprisingly, there is no revert to cartography.

Some years ago, a Google search found only one site that gave without any justification or reference, the number of ways to color the map of USA. This number was claimed to be somewhere between 20 and 21 trillion, but soon taken off the web. This left us wondering what exactly this number is. The purpose of this paper is to compute, just out of curiosity, this number for the metropolitan regions of France and the lower 48 states of the USA and prove some theorems the authors discovered on chromatic polynomials in their attempts to compute this number, for France and for the USA. 
\section{Preliminaries}
A \textit{graph} is a pair $(V,G)$ with a non-empty set $V=V(G)$ of its \textit{vertices} and a subset $E=E(G) \subset P_2(V)$ of its \textit{edges}, $E$ possibly empty. Here $P_2(V)$ is the set of all sets $T\subset V$ with the carnality $|T|=2$. To fix the notation, the following examples of graphs will appear throughout this paper. It is convenient to represent the elements of V by dots and also call them \textit{points}, whereas an edge $\{u,v\}$ in $E(G)$ will be represented by \includegraphics[width=0.9in]{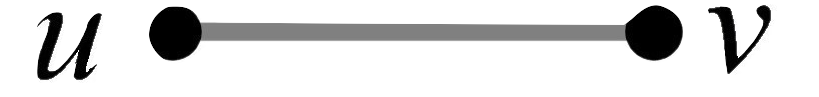}. We call \textit{u,v} \textit{adjacent} or \textit{neighbors}. All the Graphs $G$ considered in this paper are \textit{finite}, i.e. $|V(G)| < \infty$.\\\\
1. \textit{Path}. $P_n$ \includegraphics[width=1.5in]{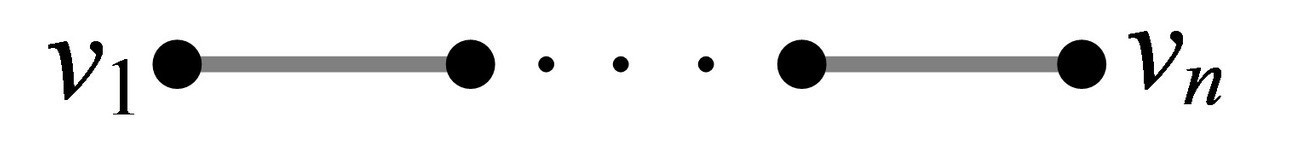} of length  $n-1 (n\geq 2)$. We call $v_1$ and $v_n$ its \textit{endpoints}.\\\\
2. \textit{Cycle}. $C_n$ on $n$ vertices $(n \geq 3)$ is obtained from $P_{n+1}$ by identifying its end points.
\begin{figure}[H]
    \centering
    \includegraphics[width=2in]{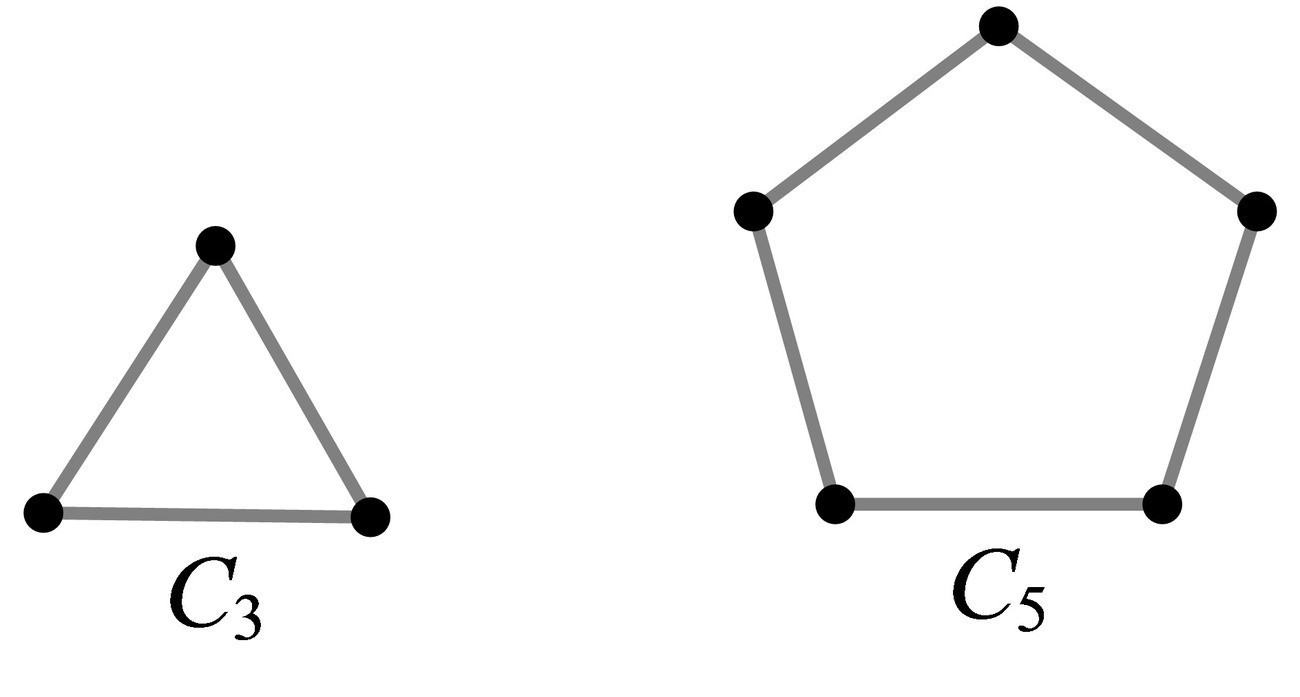}
    \caption{Cycles}
    \label{fig:my_label}
\end{figure}
\noindent
3.\hspace{2mm}\textit{Complete graph}. $G = K_n$ Here $|V(G)| = n $ and $E(G) = P_2(V)$, i.e. every two vertices $u,v$ of $G$ are adjacent.\\\\
4. \textit{Wheel}. $W_n$ on $n$ vertices $(n\geq 4)$.
\begin{figure}[H]
    \centering
    \includegraphics[width=2in]{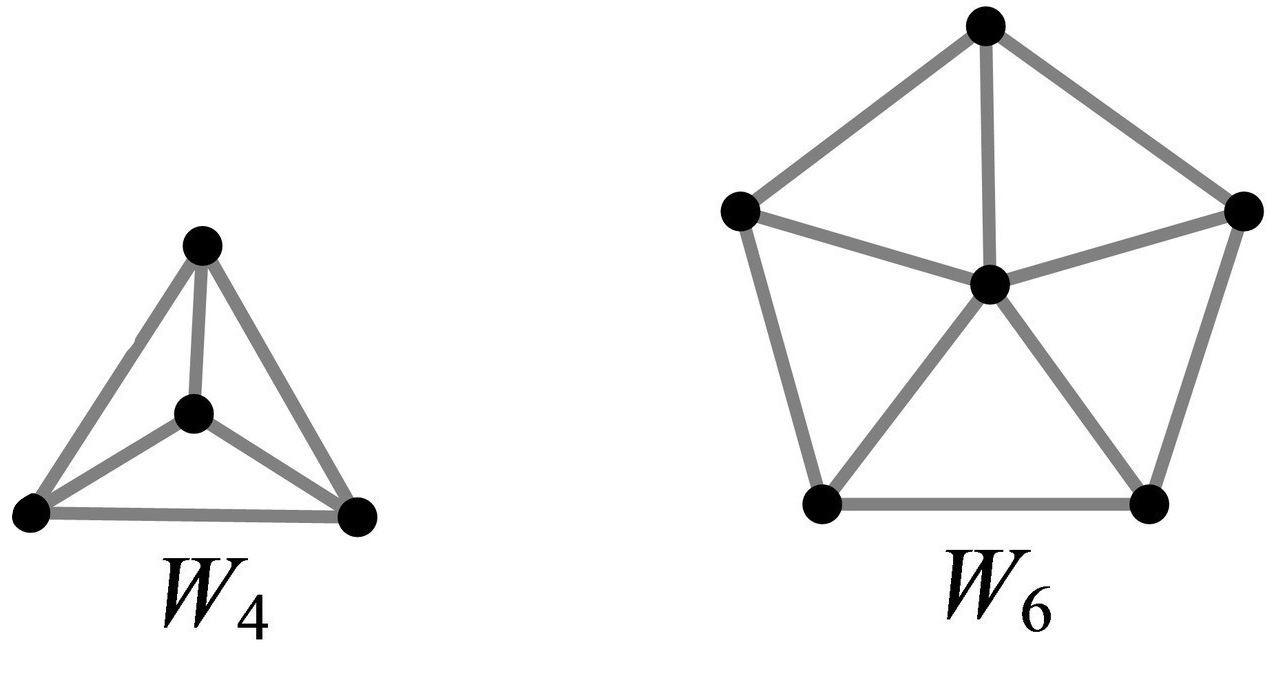}
    \caption{Wheels}
    \label{fig:my_label}
\end{figure}
\noindent
Two graphs $G_1 = (V_1, E_1)$ and $G_2 = (V_2, E_2)$ are \textit{isomorphic}, written as $G_1 \cong G_2$, if there is a bijection $f: V_1 \longrightarrow V_2$ such that $\{u,v\} \in E_1 \Longleftrightarrow \{f(u), f(v)\} \in E_2$. For example, $P_2 \cong K_2$ and $K_4 \cong W_4$. \\\\
5.\hspace{2mm}\textit{Real Life Graphs}. By these we mean the graphs associated to, among others, the geographic entities such as the map of the USA. In this graph, the vertices are the states of the USA and two states are adjacent if they share a non-trivial border, a border of positive length. Thus, for our purpose, the states of Arizona and Colorado are not adjacent, and can be given the same color.\\\\
A graph $G$ is \textit{connected} if for any two vertices $u,v$ in $V(G)$, there is a path in $G$ with $u,v$ as its end points. A \textit{subgraph} of $G = (V,E)$ is a graph $(V',E')$ with $V' \subseteq V$ and $E' \subseteq P_2(V')\cap E$. A largest connected subgraph of a graph is its \textit{connected component}. For example, the graph of the USA has 3 connected components. We will be concerned only with its largest connected components $G_A$, comprising of its lower 48 states.

\section{Coloring Problem}
\label{sec:coloringproblem}
\begin{figure}[H]
    \centering
    \includegraphics[width=3.2in]{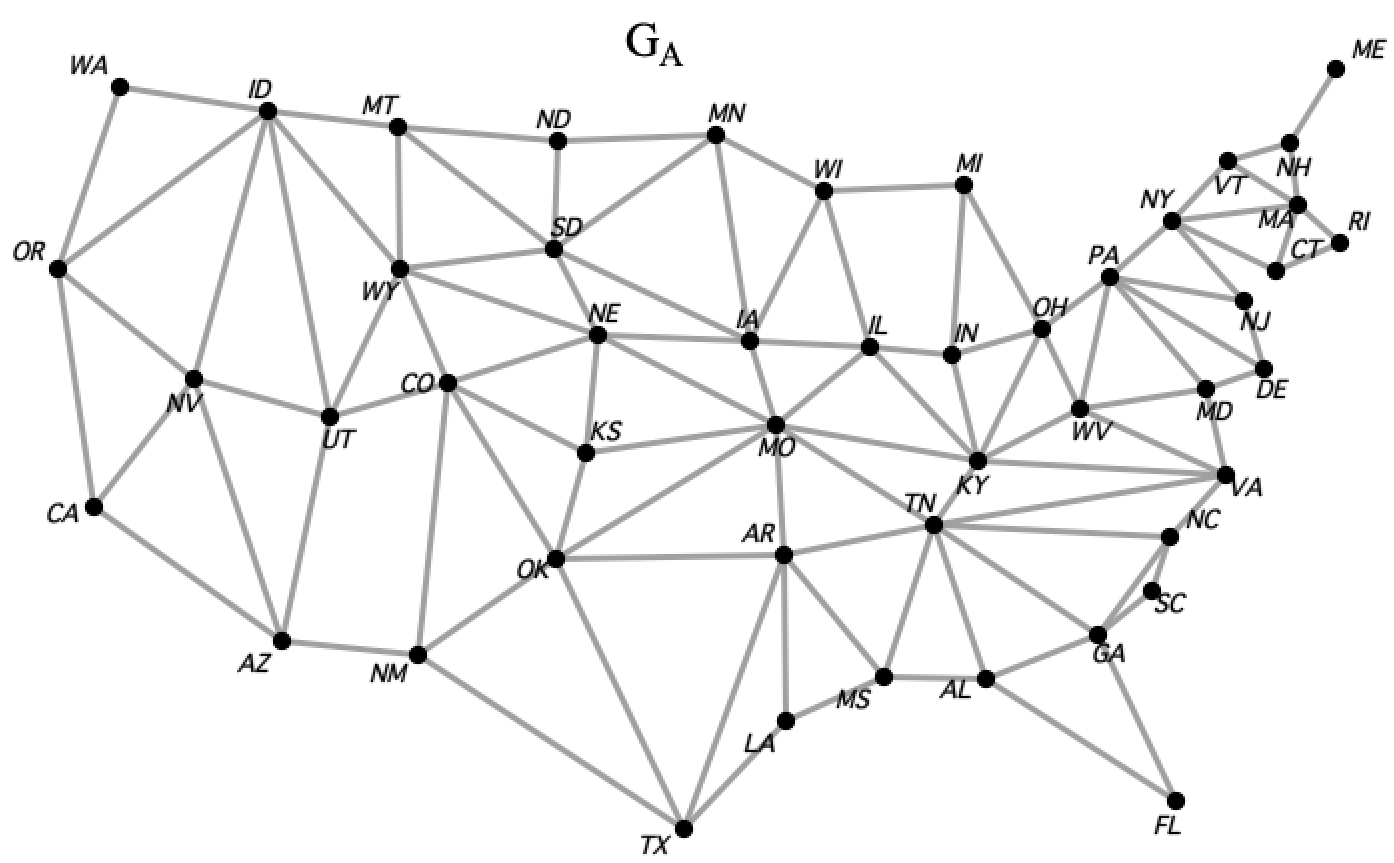}
    \caption{The graph $G_A$ with 48 vertices and 105 edges}
    \label{fig:Graph of USA}
\end{figure}

A \textit{coloring} of a graph $G$ with $t$ colors is simply a map $f: V(G) \longrightarrow C$, where $C$ is a set of given $t$ colors. The coloring is \textit{proper} if $f(u) \neq f(v)$ whenever $\{u,v\}\in E(G)$, i.e. no two adjacent vertices are given the same color. Throughout this paper, by coloring we mean only proper coloring. \\\\
The \textit{chromatic number} of a graph is the cardinality $k(G) = |C|$ of a smallest set $C$ of colors such that there is at least one way $f: V(G) \longrightarrow C$ to color it with $k(G)$ colors in $C$. The following examples will be used throughout our presentation.\\\\
1. $k(P_n) = 2$. \\\\
2. $k(C_n) = 
 \begin{cases} 
      2 & \text{if \textit{n} is even} \\
      3 & \text{if \textit{n} is odd.}
   \end{cases}$\\\\\\
3. $k(W_6) = 4.$\\\\
4. If $G_C$ is the graph of the map of Canada, it is easy to check that $k(G_C) = 3$.\\\\
5. Since the graph $G_A$ of the lower 48 states of America contains a $W_6$ centered at Nevada (cf. Figure \ref{fig:Graph of USA}) as a subgraph, $k(G_A)\geq 4$. By the Four Color Theorem, $k(G_A) = 4$.\\\\
A graph is \textit{planar} if it can be  drawn on a plane or equivalently on a sphere, without any edges crossing. All geographic maps with contiguous parts (states, provinces, districts, etc.) are planar. The Petersen graph is an example of a non-planar graph, but it does not concern us as all graphs considered here are planar. \\\\
We denote the number of ways to color a graph $G$ by a set $C$ of $t$ colors by $\chi(G,t)$. Clearly it is the cardinality of the set \textbraceleft $f:V(G) \longrightarrow C | \text{ \textit{f} is a proper coloring}$ \textbraceright. It turns out that the function $\chi(G,t)$ has the following properties(cf. \cite{Combine}, \cite{Intro}). \\\\
1. \textit{It is given by a polynomial in $t$ over $\mathbb{Z}$ (called the \textit{chromatic polynomial} of $G$).}\\\\
2. \textit{The degree $n$ of $\chi(G,t)$ is $|V(G)|$.}\\\\
3. \textit{$\chi(G,t)$ is monic, i.e. its leading coefficient is 1.}\\\\
4. \textit{The coefficient of $t^{n-1}$ is $ -|E(G)|$.}\\\\
5. \textit{The constant term of $\chi(G,t)$ is zero.}\\\\
6. \textit{The coefficients of $\chi(G,t)$ are alternately positive and negative.}\\\\
7. \textit{The degree of the lowest term in $\chi(G,t)$ is the number of components of $G$.}\\\\
8. \textit{If $|V(G)| \geq 2$, the sum of all the coefficients of $\chi(G,t)$ is zero.}\\\\
For convenience of the reader, we prove most of the properties of $\chi(G,t)$ and leave the rest as exercises.\\
\begin{proof}[Proofs]
Properties 5 and 8 are trivial. The constant term = $\chi(G,0) = 0$, because there is no way to color $G$ with no colors. This proves 5. For 8, the sum of the coefficients of $\chi(G,t)$ is the same as $\chi(G,1)$, which is zero - since there is no way to color a connected graph $G$ with 1 color if $|V(G)| \geq 2$. 
\end{proof}
\noindent
The proof of property 1, 2, 3, and 6 we will be given in Section \ref{sec:ChromReThm}. 

\section{Chromatic Polynomials of Some Standard Graphs}
\label{sec:StandardChrom}
\begin{figure}[H]
    \centering
    \includegraphics[width=3.2in]{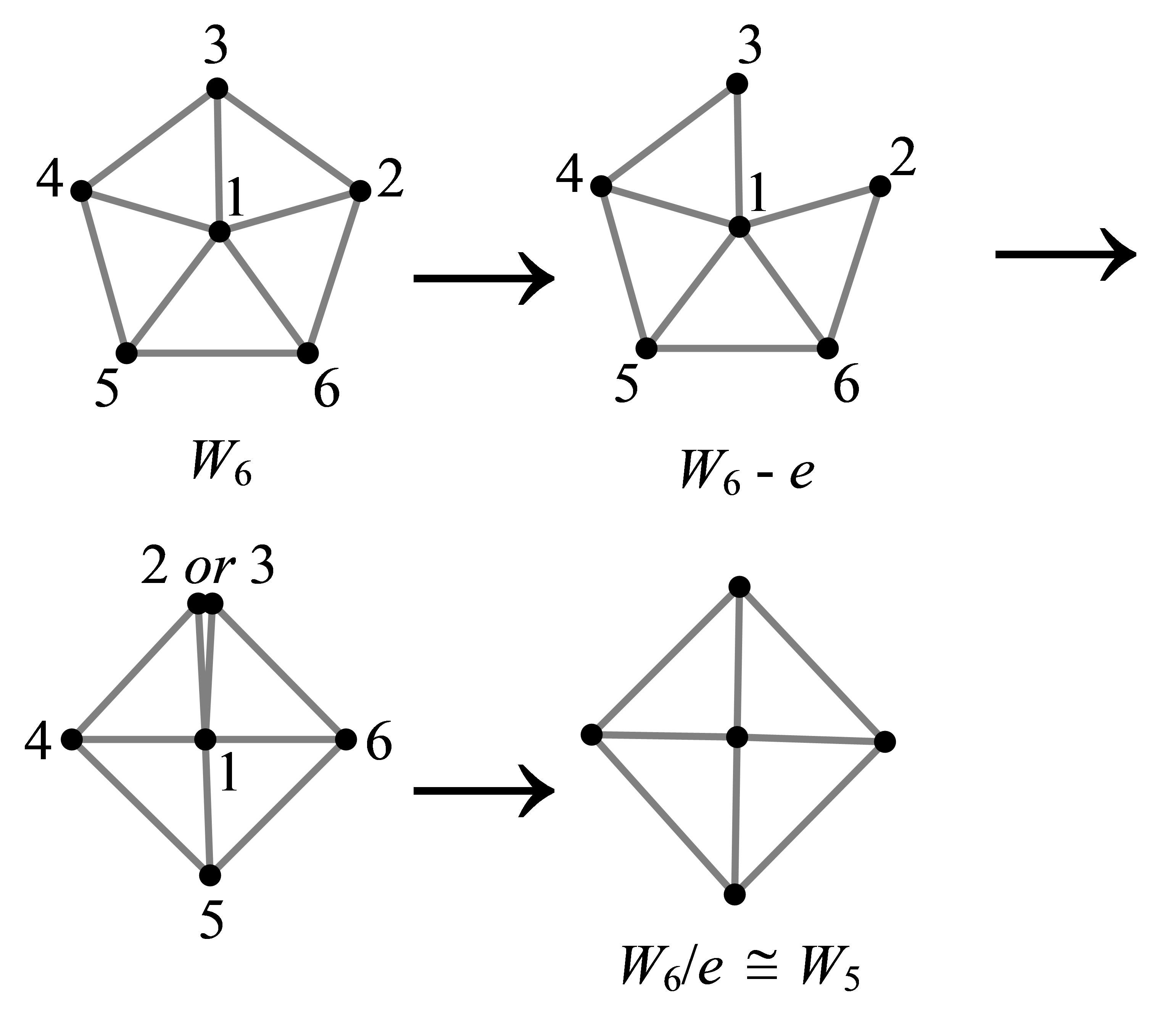}
    \caption{The broken wheel $W'_6$ and $W_6/e$}
    \label{fig:modExample}
\end{figure}

\begin{figure}[H]
    \centering
    \includegraphics[width=4.5in]{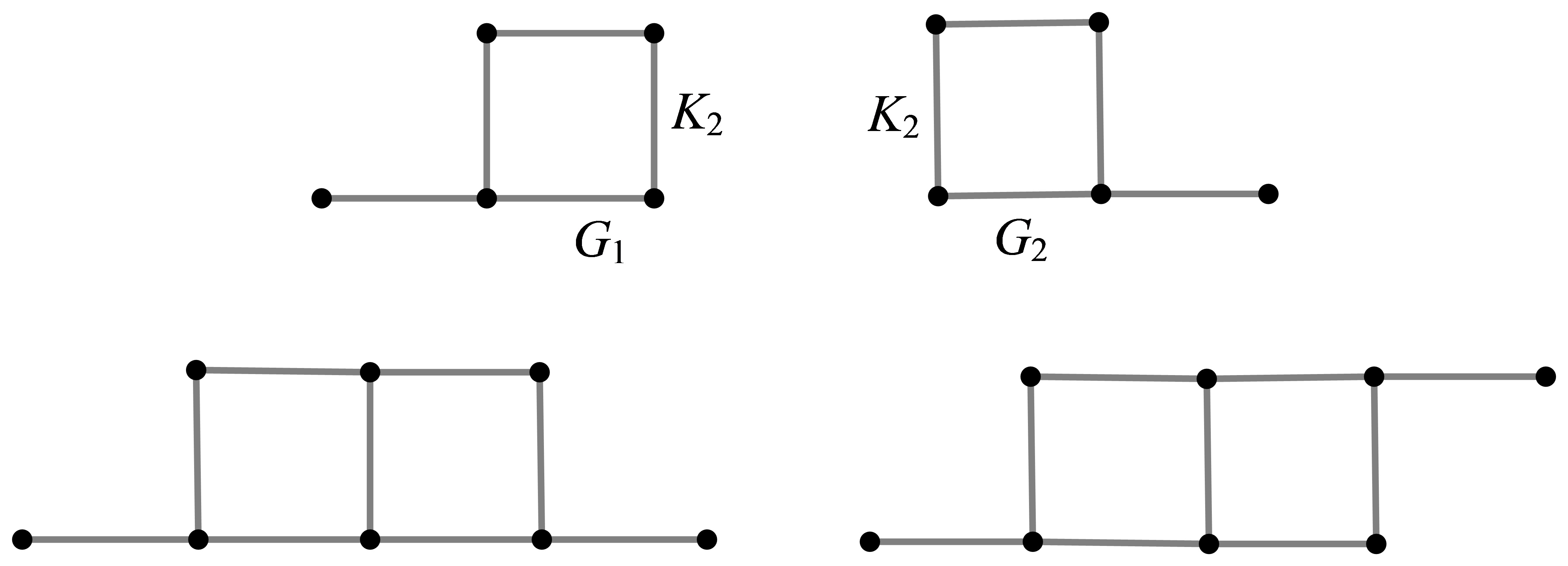}
    \caption{Non-Isomorphic overlaps}
    \label{fig:overlap}
\end{figure}

The following is obvious.\\\\
\textbf{Theorem 1}. \textit{We have}\\\\
1) $\chi(P_n,t) = t(t-1)^{n-1}$,\\
2) $\chi(C_3,t) = t(t-1)(t-2)$,\\
3) $\chi(K_n,t) = t(t-1)(t-2)...(t-n+1)$.\\\\
The following requires proofs. A connected graph $T$ is a \textit{tree} if it contains no cycle as a subgraph. \\\\
\textbf{Theorem 2}. \textit{If $T$ is a tree with $|V(T)| = n$, then $\chi(T,t) = t(t-1)^{n-1}$.}\\\\
\textbf{Theorem 3}. \textit{If $n \geq 3$, then $\chi(C_n,t) = (t-1)^{n}+(-1)^n(t-1)$.}\\\\
\textbf{Theorem 4}. \textit{If $n \geq 4,$ then}
\begin{equation}\label{eq:1}
\chi(W_n,t) = t[(t-2)^{n-1}+(-1)^{n-1} (t-2)]
\end{equation}
We will prove Theorems 3 and 4 for which, we will need the following two simple but powerful results, known as the Chromatic Reduction Theorems. In order to state and prove them, we recall some notation and terminology. If $e = \{u,v\}$ is an edge of a graph $G$, then $G-e$ is the graph $(V(G), E(G)-e)$, and $G \text{ mod } e$ is the graph $G/e$ obtained from $G-e$ by identifying $u$ and $v$. In the pictorial representation of $G/e$, we keep only one line between $u$ and $v$ if their identification creates extra lines between them. \\\\
\textbf{Example 1}. Let $e$ be the peripheral edge $\{2,3\}$ of the wheel $W_6$. Then the broken wheel $W'_6 = W_6 - e$ and $W_6/e$ are illustrated in Figure \ref{fig:modExample}. \\\\
\noindent
Now suppose that $G_1$ and $G_2$ are two graphs on disjoint sets of vertices and both $G_1$ and $G_2$ contain subgraphs that are isomorphic to a complete graph $K_l$. A graph $G$ is an \textit{overlap} of $G_1$ and $G_2$ in $K_l$ if it is obtained by identifying the subgraphs of $G_1$ and $G_2$ that are isomorphic to $K_l$.\\\\
\textbf{Example 2}. Figure~\ref{fig:overlap} illustrates two non-isomorphic overlaps of $G_1$ and $G_2$ in $K_2$.
\noindent
\textbf{Chromatic Reduction Theorem 1}. (CRT-1) If $G$ is an overlap of $G_1$ and $G_2$ in $K_l$, then
\begin{equation}\label{eq:2}
\chi(G,t) = \frac{\chi(G_1,t) \chi(G_2,t)}{\chi(K_l,t)} .
\end{equation}
To prove it, just note that the denominator in the right hand side of (\ref{eq:2}) takes care of the repeated counting in the numerator due to the overlap in $K_l$.\\\\
\textbf{Chromatic Reduction Theorem 2}. (CRT-2) If $e = \{u,v\}$ is an edge of $G$, then 
\begin{center}
$\chi(G,t) = \chi(G-e,t) - \chi(G/e,t)$.
\end{center}
To prove it, note that $\chi(G-e,t)$ is the number of ways to color $G$ when \textit{u,v} get the same color, which is $\chi(G/e,t)$, plus the number of ways to color $G$ when $v,u$ get different colors, which is $\chi(G,t)$.\\\\
\textbf{Remarks}. \\\\
1. It is easy to construct examples to show that CRT-1 is false if the overlap is not a complete graph. \\\\
2. Example 2 shows that two non-isomorphic graphs can have the same chromatic polynomial.\\\\
3. To compute the chromatic polynomial $\chi(G,t)$, there is no loss of generality to assume that $G$ is connected, because if $G_1$ and $G_2$ are its two connected components then $\chi(G,t)=\chi(G_1,t) \chi(G_2,t)$. 

\section{Application of Chromatic Reduction Theorems}
\label{sec:ChromReThm}
As an easy application of the Chromatic Reduction Theorems, first we will prove Theorems 3 and 4. The proofs of Theorems 1 and 2 are easy and are left as an exercise for the reader.

\begin{proof}[Proof of Theorem 3]
We use induction on $n$. If $n = 3$, the right hand side of the equation
\begin{center}
  $\chi(C_n, t) = (t-1)^n +(-1)^n (t-1)$
\end{center}
is $(t-1)^3 - (t-1) = (t-1)((t-1)^2 - 1) = t(t-1)(t-2)$ which is $\chi(C_3, t)$. If $n\geq3$, let $e$ be an edge of $C_{n+1}$. Then $C_{n+1} - e \cong P_{n+1}$ and $C_{n+1}/e \cong C_n$.\\
By CRT-2 and the Induction Hypothesis,
\begin{align*}
 \chi(C_{n+1}, t) 
    & = \chi(P_{n+1},t) - \chi(C_n,t)\\
    & = t(t-n)^n - [(t-1)^n + (-1)^n (t-1)]\\
    & = (t-1)^{n+1} + (-1)^{n+1} (t-1)
\end{align*}
This completes the proof by induction.
\end{proof}
\noindent
In order to prove Theorem 4, we use the following terminology (See Figure \ref{fig:modExample}). For $n \geq 5$, let $e$ be an edge on the rim of $W_n$. We call $W'_n = W_n - e$ a \textit{broken wheel}. Clearly, $W_n / e \cong W_{n-1}$. We also need the following.\\\\
\textbf{Lemma}. \textit{For $n\geq 4$, $\chi(W'_{n+1}, t)= t(t-1)(t-2)^{n-1}$}.
\medskip
\begin{proof}[Proof of Lemma]
Clearly, $W'_{n+1}$ is an overlap of $C_3$ and $W'_n$ in $K_2$. By CRT-1 and the Induction Hypothesis, 

\begin{align*}
    \chi(W'_{n+1}, t)
    & = \frac{\chi(W'_n,t) \chi(C_3,t)}{\chi(K_2,t)}\\
    & = \frac{(t (t-1) (t-2)^{n-2})(t (t-1) (t-2))}{t (t-1)}\\
    & = t(t-1)(t-2)^{n-1}.
\end{align*}
\end{proof}

\begin{proof}[Proof of Theorem 4]
 It is again by induction on $n$. First let $n=4$. Since $W_4 \cong K_4$, $\chi(W_4,t) = \chi(K_4,t) = t(t-1)(t-2)(t-3)$, which is the same as the right hand side of the equation (1) for $n=4$. Now let $n\geq 4$. Then by CRT-2 and the Induction Hypothesis,
 
 \begin{align*}
     \chi(W_{n+1},t)
     & = \chi(W'_{n+1},t) - \chi(W_n,t)\\
     & = t(t-1)(t-2)^{n-1} - t[(t-2)^{n-1} + (-1)^{n-1}(t-2)]\\
     & = t[(t-2)^{n} + (-1)^n (t-2)].
 \end{align*}
\end{proof}
\noindent
We now combine some of the properties of the chromatic polynomial as the following result.\\\\
\textbf{Theorem 5}. \textit{Let $G$ be a connected graph on $n$ vertices. The counting function $\chi(G,t)$ is given by a monic polynomial over $\mathbb{Z}$ of degree $n$. If $n \geq 2$, its coefficients are alternately positive and negative and their sum us zero.}
\medskip
\begin{proof}
 We prove it by double induction on $n = |V(G)|$ and $k = |E(G)|$. It is easy to check that it is true for all graphs with $ n \leq 3$. Let $n \geq 3$ be any integer such that it is true for all graphs with $n \geq 3$. We will prove it for all graphs with $|V(G)| = n+1$. It is certainly true if $k = 0$. We now apply induction to $ k = |E(G)|$. Suppose it is true for all graphs with $|V(G)| = n+1$ and $|E(G)| \leq k$ for a given $ k\geq 0$. If $G$ is a graph with $|V(G)| = n+1$ and $|E(G)| = k+1$, let $e \in E(G)$. By the Induction Hypothesis, 
 \begin{align*}
     \chi(G-e, t) 
     & = t^{n+1} - a_{n}t^{n} + a_{n-1}t^{n-1} - ... \hspace{5mm},\\
     and \hspace{5mm} \chi(G/e, t) 
     & = t^{n} - b_{n-1}t^{n-1} + ...
    \text{\hspace{10mm}with } a_j, b_j \text{ all positive integers.} 
 \end{align*}
\noindent
 Therefore by CRT-2,
 \begin{align*}
     \chi(G,t)
     & = \chi(G-e,t) - \chi(G/e,t)\\
     & = t^{n+1} - (a_{n} + 1)t^{n} + (a_{n-1}+ b_{n-1})t^{n-1} - ...
 \end{align*}
 This completes the proof by double induction.
\end{proof}
\noindent
 That the sum of the coefficients is zero, has already been proved earlier in Section \ref{sec:coloringproblem}.
 
\section{Chromatic Polynomials of Geographic Graphs}
\begin{figure}[H]
    \centering
    \includegraphics[width=4in]{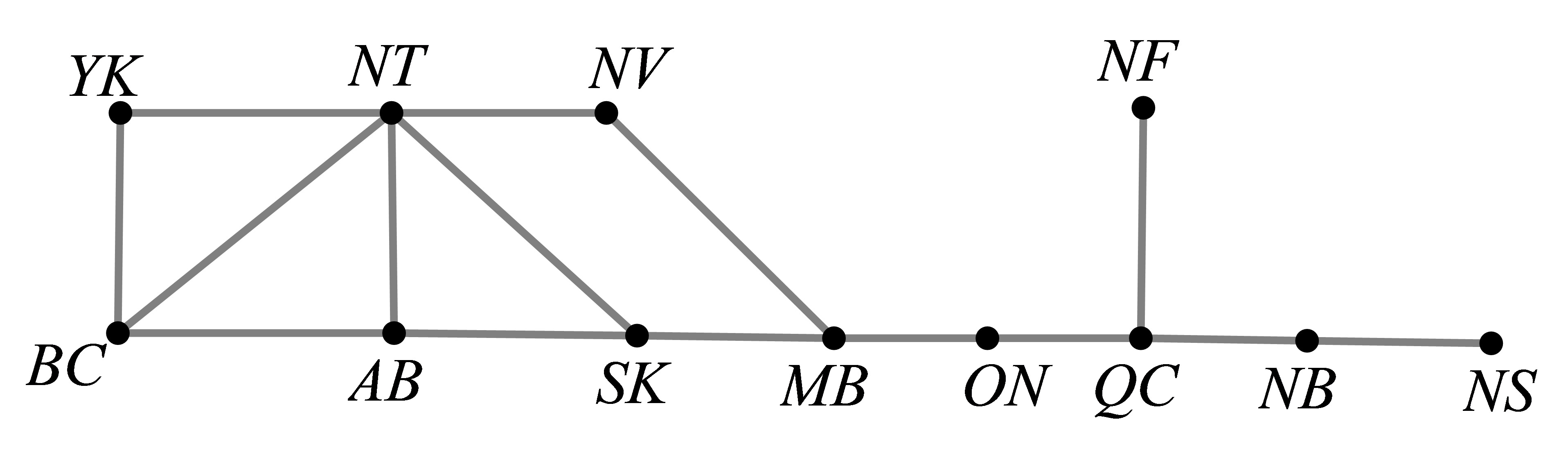}
    \caption{The graph $G_C$ of the map of Canada}
    \label{fig:Canada}
\end{figure}

\begin{figure}[H]
    \centering
    \includegraphics[width=2in]{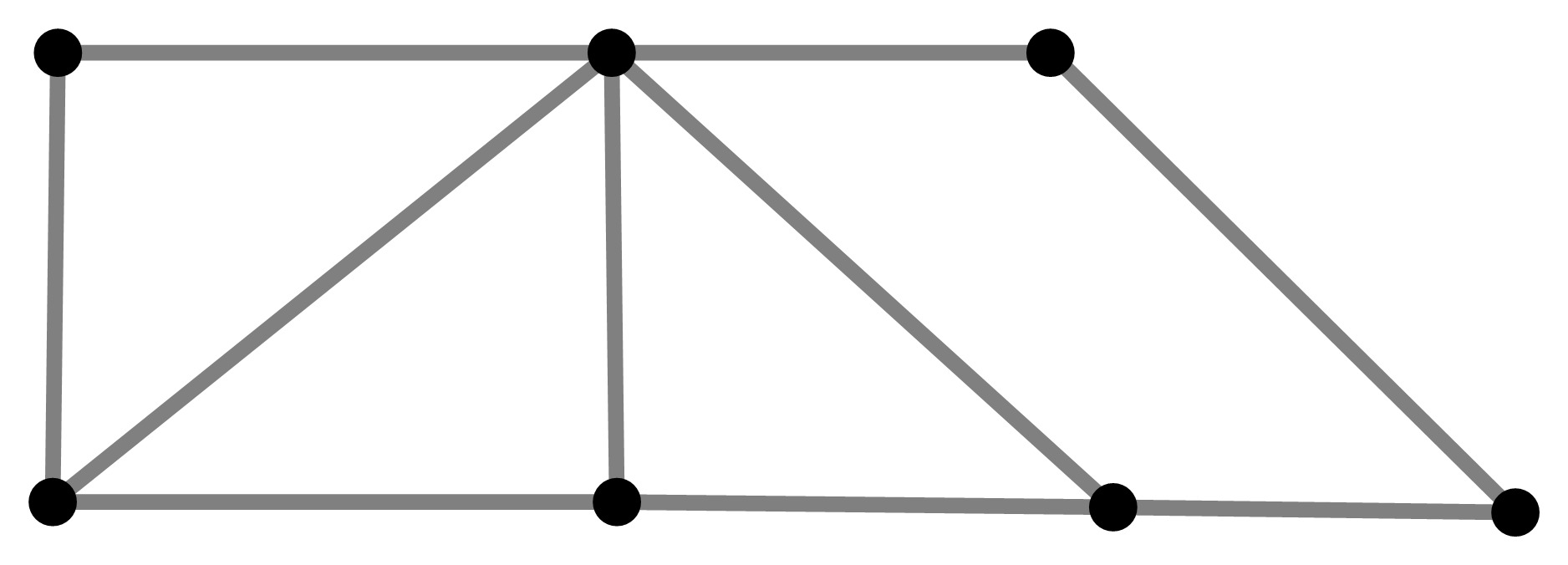}
    \caption{The graph of K}
    \label{fig:graphK}
\end{figure}

In this section, we will compute the chromatic polynomials of the graphs that are associated with the geographic maps of Canada, France, and America.\\\\
1.\hspace{2mm}\textbf{Canada}.
    First we take the easy and almost trivial case of the map of the provinces and territories of
    Canada.\footnote{Map of Canada retrieved from: https://www.conceptdraw.com/How-To-Guide/geo-map-canada-prince-edward-island} Its graph $G_C$ is illustrated in Figure \ref{fig:Canada}, 
\begin{figure}[H]
    \centering
    \includegraphics[width=4in]{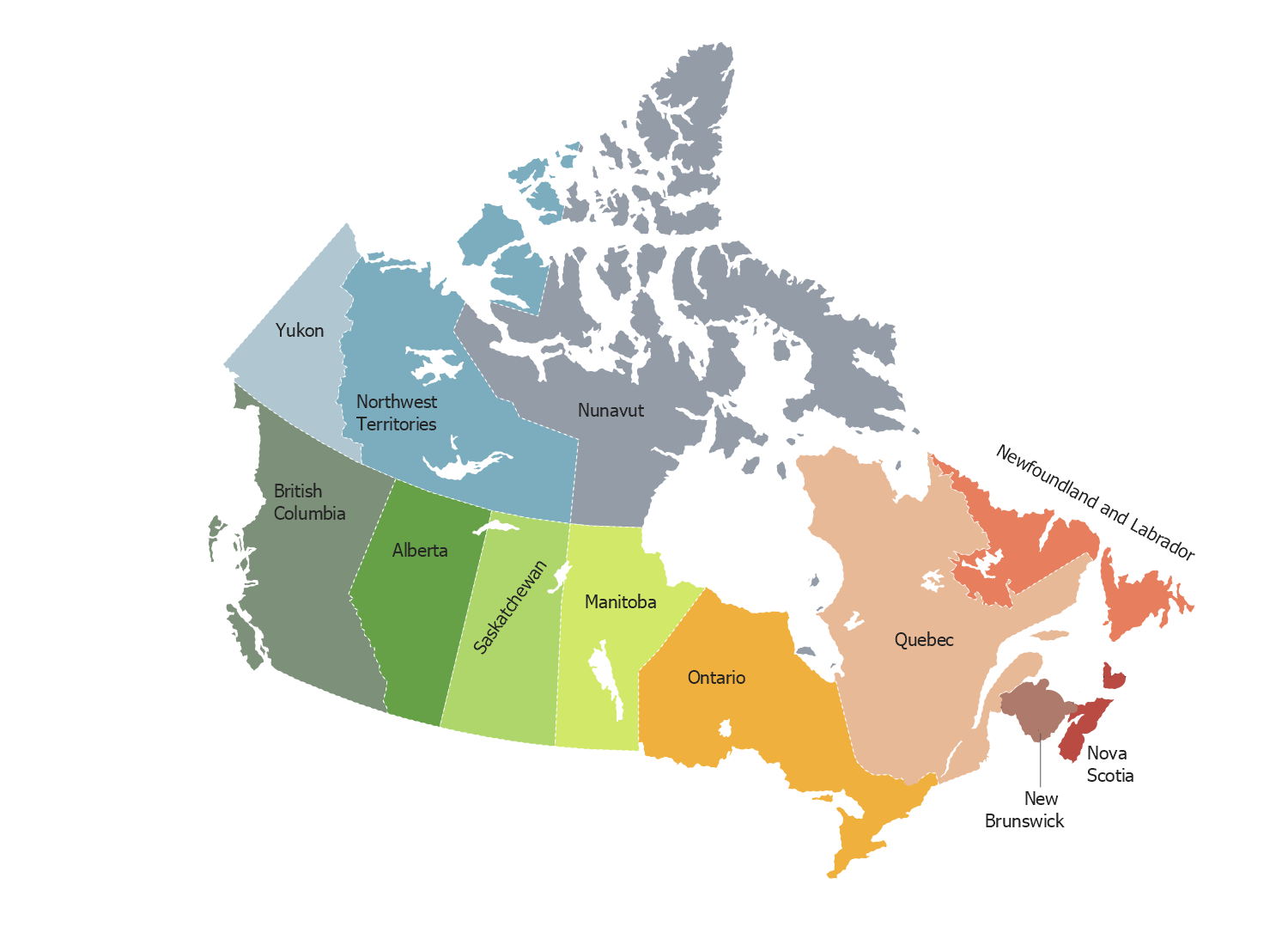}
    \caption{Map of Canada}
    \label{fig:map_Canada}
\end{figure}
\noindent
which is an overlap of a tree $T$, of 6 vertices and the graph $K$ (see Figure~\ref{fig:graphK}) in $K_1$ (at the province of Manitoba).
\noindent
The graph $K$ is an overlap of $C_4$ in $K_2$ and a subgraph which is a series of repeated overlaps of $C_3$ in $K_2$. By CRT-1, Theorem 2 and Theorem 3, $\chi(G_C,t) = \frac{\chi(T,t)\chi(K,t)}{\chi(K_1,t)}$. But $\chi(K,t) = \frac{\chi(C_4,t){\chi(C_3,t)}^3}{{\chi(K_2, t)}^3}$. Therefore, $\chi(G_C,t) = t(t-1)^6(t-2)^3(t-3t+3)$.
Thus the \\\\
\fbox{
\parbox{\textwidth}{
\textbf{number of ways to color the map of Canada with 3 colors is $\boldsymbol{\chi(G_C,3) = 576}.$}
}
}
\\\\

\begin{figure}[H]
  \centering
  \includegraphics[width=3in]{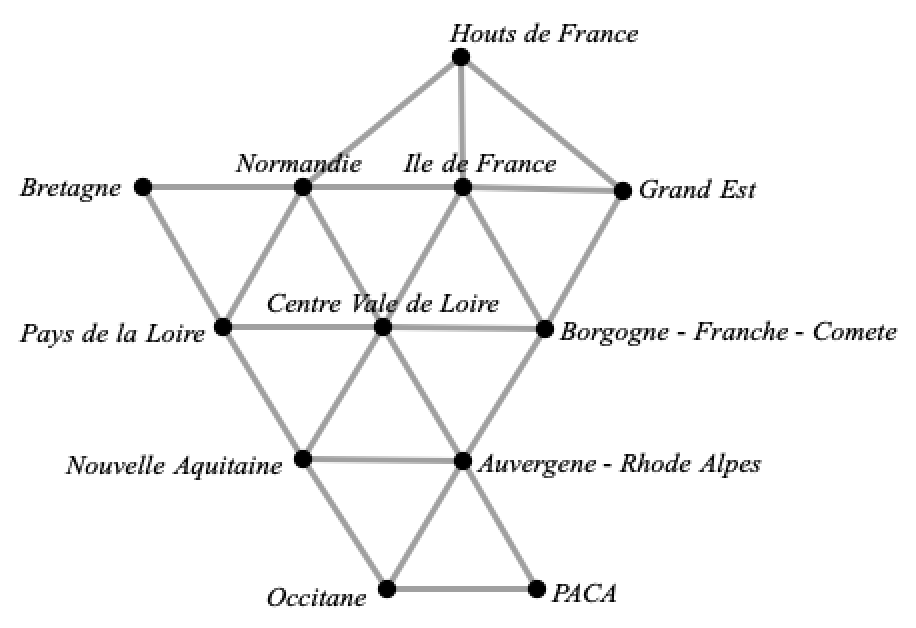}
  \captionof{figure}{The graph $G_F$ of the 12 contiguous regions of France}
  \label{fig:France}
\end{figure}

\noindent
2. \textbf{France}.
The country of France has 18 regions, but we consider only the 12 contiguous regions.\footnote{Map of France retrieved from: http://evasion-online.com/tag/carte-des-regions-de-france-2017}, to which we associate its graph $G_F$, as illustrated in Figure~\ref{fig:France}.
\begin{figure}[H]
    \centering
    \includegraphics[width=2in]{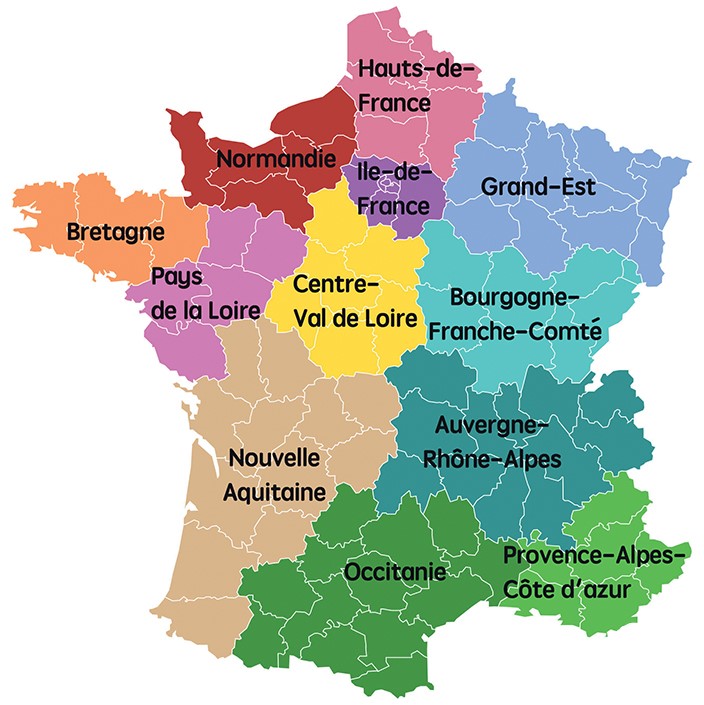}
    \caption{Map $G_F$ of France}
    \label{fig:map_France}
\end{figure}

\begin{figure}[H]
    \centering
    \includegraphics[width=2in]{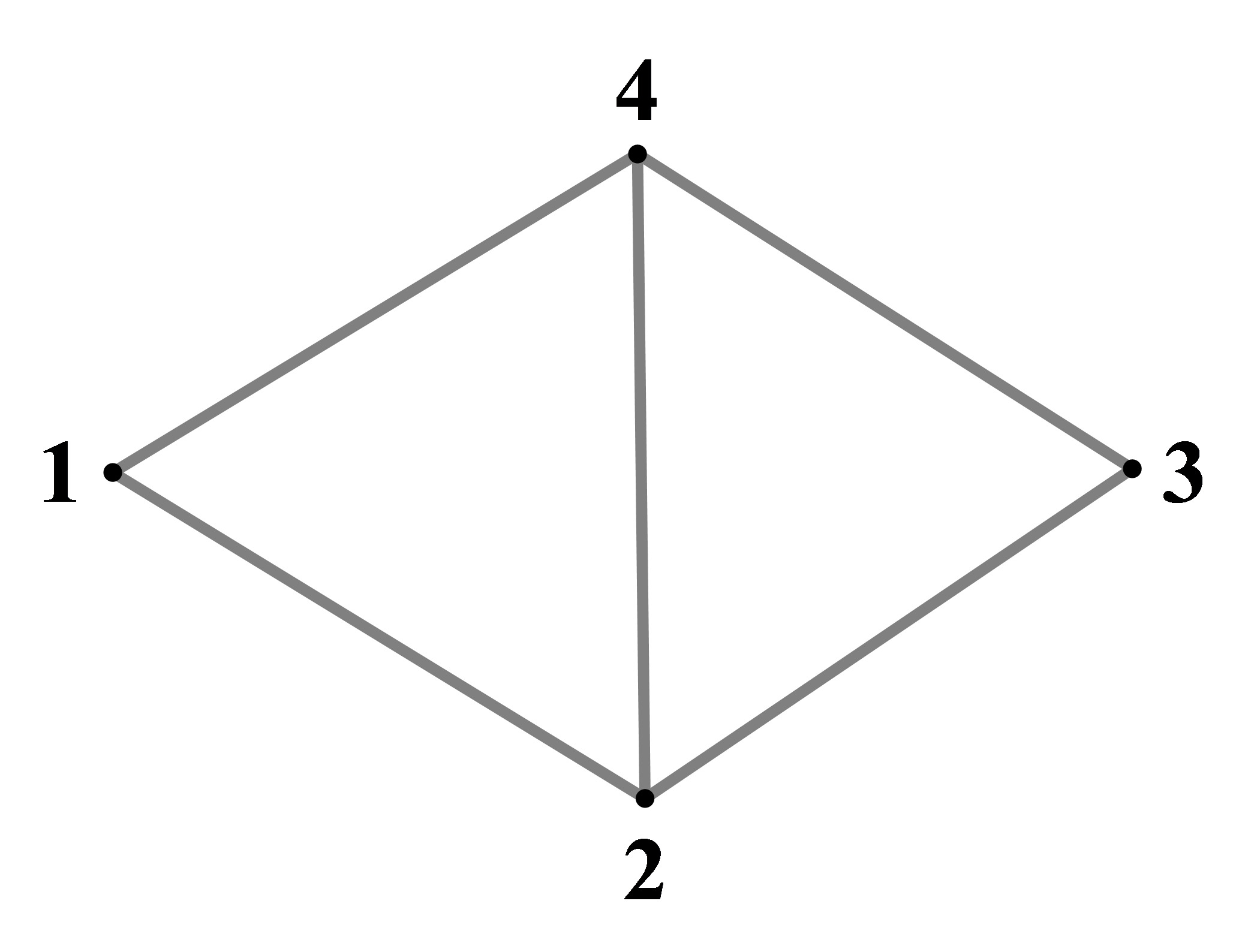}
    \caption{Double wedge $\wedge_2$}
    \label{fig:Double}
\end{figure}

\begin{figure}[H]
    \centering
    \includegraphics[width=2in]{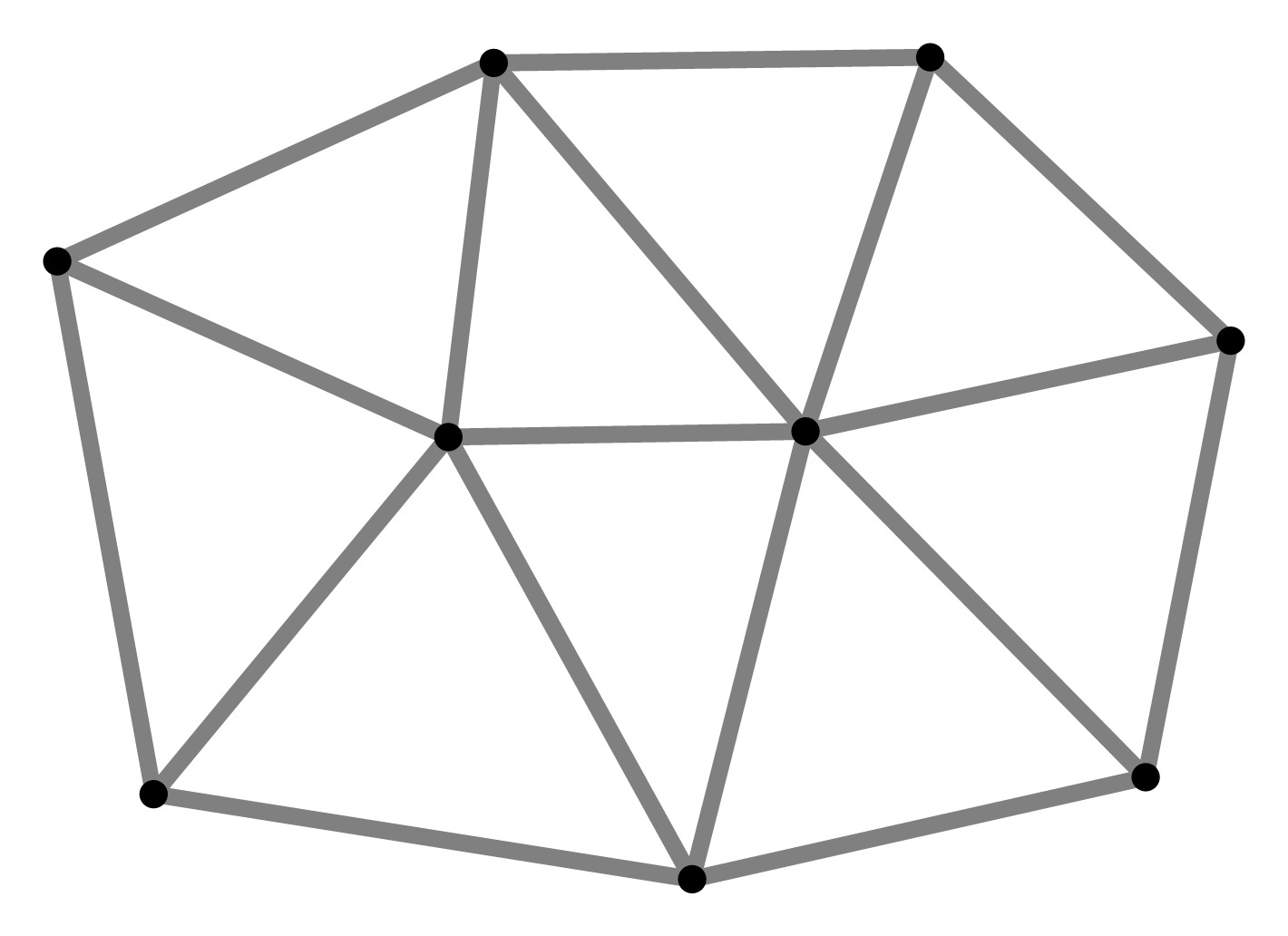}
    \caption{The interlocking wheel $W_6\wedge_2W_7$}
    \label{fig:interlock}
\end{figure}

The graph $G_F$ contains not only cycles, but also wheels. Computing $\chi(G_F, t)$ would have been equally easy if they overlapped in complete graphs. Instead the wheels "interlock" in a double wedge $\wedge_2$ (Figure \ref{fig:Double} and \ref{fig:interlock}). The computer programs delete and mod out edges and combine the chromatic polynomials of the loads of graphs appearing in the process (cf. Figure \ref{fig:slide}) but it is exceedingly difficult to do it by hand. We tried in vain, but ultimately had to use a computer program to put the pieces together. 
\begin{figure}[H]
    \centering
    \includegraphics[width=6in]{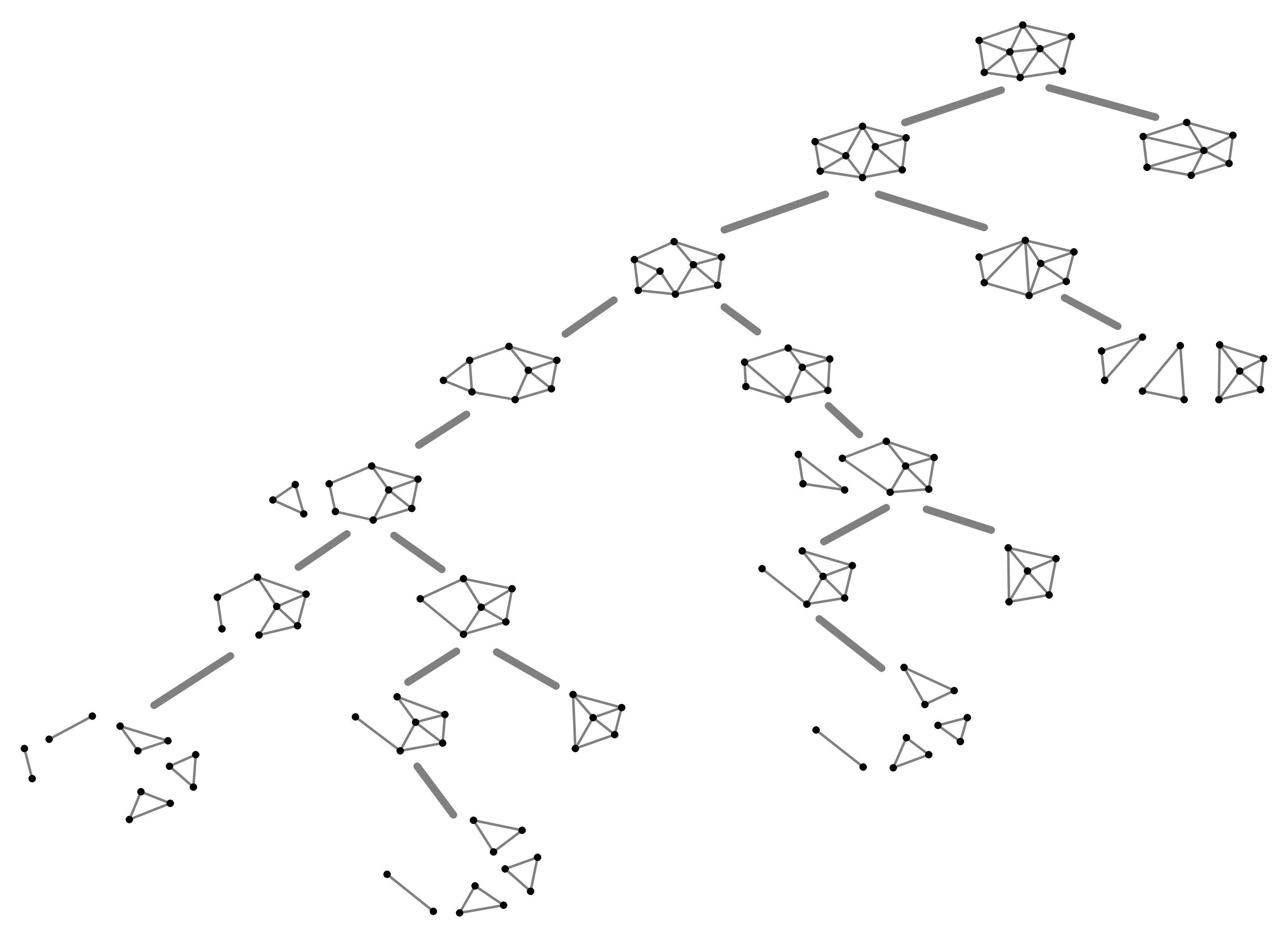}
    \caption{Delete and mod out surgery on $W_6 \wedge_2 W_6$}
    \label{fig:slide}
\end{figure}
\noindent
One could add the edge $e = \{1,3\}$ to make it complete and apply CRT-2, but when we mod out edges, in the process we end up with graphs that are in general not even planar. Thus computing $\chi(G_F,t)$ requires some new ideas.

To overcome this difficulty, we prove the following result on interlocking wheels, which resulted in an exponential simplification in our calculations. First we explain our notation. We denote by $W_m\wedge_2 W_n$ the graph which is an overlap of two wheels $W_m$ and $W_n$ in the wedge $\wedge_2$. See Figure \ref{fig:interlock} for $m=6$ and $n=7$.

\noindent
\textbf{Main Theorem}. If $m,n \geq 4$ and $m+n \geq 9$, then

\begin{equation}\label{eq:3}
\begin{split}
\chi(W_m \wedge W_n,t) 
& = \frac{1}{t-1}[t(t-2)(t-3)((t-2)^{n-3} + (-1)^n)((t-2)^{m-3} + (-1)^m) \\
&  + t (t-2)^3((t-2)^{n-4} + (-1)^{n-1})((t-2)^{m-4} + (-1)^{m-1}]
\end{split}
\end{equation}

\begin{proof} 
First we remark that by theorem 5, the polynomial in the bracket must have $t-1$ as a factor.

Our proof is by two way induction on $n$ and $m$.
\newline Base Case:
For $n = 5$ and $m = 4$, by CRT-2,
\begin{equation*}
\begin{split}
\chi(W_5 \wedge_2 W_4,t)&=\chi((W_5 \wedge_2 W_4)-e,t)-\chi((W_5 \wedge_2 W_4)/e,t)\\
\end{split}
\end{equation*}
But, 
\begin{equation*}
\begin{split}
\chi((W_5 \wedge_2 W_4)-e,t)&=\chi(W_5,t)=t((t-2)^4-(-1)^5(t-2))\\
\end{split}
\end{equation*}
and
\begin{equation*}
\begin{split}
 \chi((W_5 \wedge_2 W_4)/e,t) &=\frac{\chi(C_3,t)^2}{\chi(K_2,t)}\\
 &=\frac{(t(t-1)(t-2))^2}{t(t-1)}\\
 &=\frac{t^2(t-1)^2(t-2)^2}{t(t-1)}\\
 &=t(t-1)(t-2)^2\\
 \end{split}
 \end{equation*}
 Thus
 \begin{equation*}
\begin{split}
 \chi(W_5 \wedge_2 W_4,t)&=t((t-2)^4-(-1)^5(t-2)-t(t-1)(t-2)^2\\
&=t(t-1)(t-2)(t-3)^2\\ 
\end{split}
\end{equation*}
On the other hand, for $n=5$ and $m=4$. The right hand side of equation (\ref{eq:3}) is
\begin{equation*}
\begin{split}
 &=
\frac{t(t-2)(t-3)((t-2)^{2}-(-1)^{4})((t-2)-(-1)^{3})}{(t-1)}\\
&+\frac{t(t-2)^3((t-2)-(-1)^{3})((t-2)^{0}-(-1)^{2})}{(t-1)}\\
&=\frac{t(t-2)(t-3)((t-2)^{2}-(-1)^{4})((t-2)-(-1)^{3})}{(t-1)}\\
&=\frac{t(t-2)(t-3)((t-2)^{2}-(-1)^{4})(t-1)}{(t-1)}\\
&=t(t-2)(t-3)((t-2)^{2}-(-1)^{4})=t(t-1)(t-2)(t-3)^2 .\\
\end{split}
\end{equation*}
Therefore, the theorem holds for $m=4$ and $n=5$.
We check similarly that the theorem holds for $n=5$ and $m=5$. \\
By CRT-2,
\begin{equation*}
    \begin{split}
\chi(W_5 \wedge_2 W_5,t)&=\chi((W_5 \wedge_2 W_5)+e,t)+\chi((W_5 \wedge_2 W_5)/e,t)\\
\end{split}
\end{equation*}
But (see Figure \ref{fig:slide}),
\begin{figure}[H]
    \centering
    \includegraphics[width=2in]{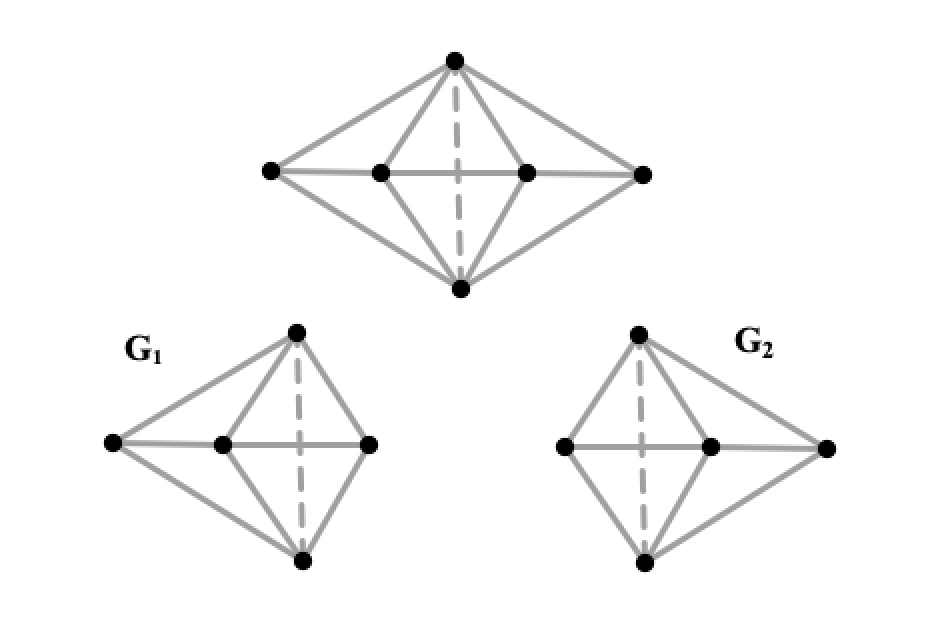}
    \caption{Graph of $G_1$ and $G_2$}
    \label{fig:slide}
\end{figure}
\begin{equation*}
\chi((W_5 \wedge_2 W_5)+e,t) =\frac{\chi(G_1,t)\chi(G_2,t)}{\chi(K_4,t)} .
\end{equation*}
But
\begin{align*}
    \chi(G_1,t)=\chi(G_2,t)&=\frac{\chi(W_4,t)\chi(C_3,t)}{\chi(K_2,3)}-\chi(W_4,t)\hspace{1mm}\\
&=\frac{t^2(t-1)^2(t-2)^2(t-3)}{t(t-1)}- t(t-1)(t-2)(t-3)\\
&=t(t-1)(t-2)^2(t-3)-t(t-1)-t(t-1)(t-2)(t-3)\\
&=t(t-1)(t-2)(t-3)^2 .\\
\end{align*}
Therefore,
\begin{align*}
\chi((W_5 \wedge_2 W_5)+e,t)&= \frac{(t(t-1)(t-2)(t-3)^2)^2}{t(t-1)(t-2)(t-3)}\\
&= \frac{t^2(t-1)^2(t-2)^2(t-3)^4}{t(t-1)(t-2)(t-3)}\\
&=t(t-1)(t-2)(t-3)^3 ,\\
\end{align*}
and
\begin{align*}
\chi((W_5 \wedge_2 W_5)/e,t)&=\frac{\chi(C_3,t)^3}{\chi(K_2,t)^2}=\frac{(t(t-1)(t-2))^3}{(t(t-1))^2}\\
&=\frac{t^3(t-1)^3(t-2)^3}{t^2(t-1)^2}\\
&=t(t-1)(t-2)^3 .\\
\end{align*}
Therefore,
\begin{align*}
\chi(W_5 \wedge_2 W_5,t)&=t(t-1)(t-2)(t-3)^3+t(t-1)(t-2)^3\\
&=t(t-1)(t-2)(t^3-8t^2+23t-23) .\\
\end{align*}
On the other hand, for $n=5$ and $m=5$, the right hand side of equation (\ref{eq:3}) is
\begin{align*}
=&\frac{t(t-2)(t-3)((t-2)^{2}-(-1)^{4})((t-2)^{2}-(-1)^{4})}{(t-1)}\\
&+\frac{t(t-2)^3((t-2)-(-1)^{3})((t-2)-(-1)^{3})}{(t-1)}\\
=&\frac{t(t-2)(t-3)({(t-2)^{2}-(-1)^{4}})^2}{(t-1)}\\
&+\frac{t(t-2)^3((t-2)-(-1)^{3})^2}{(t-1)}\\
=&\frac{t(t-2)(t-3)(t-1)^2(t-3)^2+t(t-2)^3((t-1)^2}{(t-1)}\\
=&t(t-2)(t-3)^3(t-1)+t(t-2)^3(t-1)\\
=&t(t-1)(t-2)((t-3)^3+(t-2)^2)\\
=&t(t-1)(t-2)(t^3-8t^2+23t-23).\\
\end{align*}
Thus the theorem holds for $n=5$ and $m=5$.\\\\
\noindent
Now assume that (\ref{eq:3}) is true for $n = 5$ and $m \geq 5$. We must show it's true for $n = 5$ and $m + 1 \geq 5$. \\
Again by CRT-2 and by the Inductive Hypothesis,

  \begin{align*}
  \chi(W_5 \wedge_2 W_{m+1},t)&=\frac{\chi(W_5,t)\chi(C_3,t)^{(m+1)-4}}{\chi(K_2,t)^{(m+1)-4}}-\chi(W_5\wedge_2W_m,t)\\
&=\frac{t((t-2)^4-(-1)^5(t-2))t^{m-3}(t-1)^{m-3}(t-2)^{m-3}}{t^{m-3}(t-1)^{m-3}}\\
&-\frac{t(t-2)(t-3)((t-2)^{2}-(-1)^{4})((t-2)^{m-3}-(-1)^{m-1})}{(t-1)}\\
&-\frac{t(t-2)^3(t-1)((t-2)^{m-4}-(-1)^{m-4})}{(t-1)}\\
&=t(t-2)^{m-2}((t-2)^3+1)\\
&-t(t-2)(t-3)^2((t-2)^{m-3}-(-1)^{m-1})\\
&-t(t-2)^3((t-2)^{m-4}-(-1)^{m-2})\\
&=t(t-1)(t-2)^{m-2}(t^2-5t+7)\\
&-t(t-2)(t-3)^2((t-2)^{m-3}-(-1)^{m-1})\\
&-t(t-2)^3((t-2)^{m-4}-(-1)^{m-2})\\
&=t(t-2)[(t-1)(t-2)^{m-3}(t^2-5t+7)\\
&-(t-3)^2((t-2)^{m-3}-(-1)^{m-1})\\
&-(t-2)^2((t-2)^{m-4}-(-1)^{m-2})\\
&= e_1, \text{say}.
 \end{align*}
 We want to show that if the right hand side of (\ref{eq:3}) for $n = 5$ and $m + 1$, which is 
 \begin{align*}
       & t(t-2)[(t-3)^2((t-2)^{(m+1)-3}-(-1)^{(m+1)-1})\\
       & +(t-2)^2((t-2)^{(m+1)-4}-(-1)^{(m+1)-2})]\\ 
       & = e_2 , \text{say}
 \end{align*}
then \[e_1=e_2.\]
By Wolfram Mathematica \[e_1-e_2=0\]
which implies \[e_1=e_2\].
\\~\\
Now assume (\ref{eq:3}) is true for $n\geq5$ and $m\geq5$. We must prove it true for $n+1\geq5$ and $m\geq5$. \\
First,
 \[\chi(W_{n+1} \wedge_2 W_{m},t)=\frac{\chi(W_{n+1},t)\chi(C_3,t)^{(m+1)-4}}{\chi(K_2,t)^{(m+1)-4}}-\chi(W_n\wedge_2W_m,t),\]
 which by the Inductive Hypothesis is
 \begin{align*}
     =&\frac{t(t-2)((t-2)^{(n+1)-2}-(-1)^{(n+1)-1}(t(t-1)(t-2))^{(m+1)-4}}{(t(t-1))^{(m+1)-4}}\\
     &-\frac{t(t-2)(t-3)((t-2)^{n-3}-(-1)^{n-1})((t-2)^{m-3}-(-1)^{m-1})}{(t-1)}\\
     &-\frac{t(t-2)^3((t-2)^{n-4}-(-1)^{n-2})((t-2)^{m-4}-(-1)^{m-2})}{(t-1)}\\
 =&t(t-2)^{(m-2)}((t-2)^{n-1}-(-1)^{n})\\
 &-\frac{t(t-2)(t-3)((t-2)^{n-3}-(-1)^{n-1})((t-2)^{m-3}-(-1)^{m-1})}{(t-1)}\\
 &-\frac{t(t-2)^3((t-2)^{n-4}-(-1)^{n-2})((t-2)^{m-4}-(-1)^{m-2})}{(t-1)}.\\
 \end{align*}
 Now we want to show that if 
 \begin{align*}
        e_1=&t(t-2)^{(m-2)}((t-2)^{n-1}-(-1)^{n})\\
        &-\frac{t(t-2)(t-3)((t-2)^{n-3}-(-1)^{n-1})((t-2)^{m-3}-(-1)^{m-1})}{(t-1)}\\
        &-\frac{t(t-2)^3((t-2)^{n-4}-(-1)^{n-2})((t-2)^{m-4}-(-1)^{m-2})}{(t-1)}\\  
 \end{align*}
 and
 \begin{align*}
      e_2=&\frac{t(t-2)(t-3)((t-2)^{(n+1)-3}-(-1)^{(n+1)-1})((t-2)^{(m+1)-3}-(-1)^{(m+1)-1})}{(t-1)}\\
      &+\frac{t(t-2)^3((t-2)^{(n+1)-4}-(-1)^{(n+1)-2})((t-2)^{(m+1)-4}-(-1)^{(m+1)-2})}{(t-1)}\\ 
 \end{align*}
\noindent
 then \[e_1=e_2\]
 \noindent
 Again by Wolfram Mathematica 
 \[e_1-e_2=0,\] 
 which implies \[e_1=e_2\].
 This proves (\ref{eq:3}) for $n + 1 \geq 5$ and $m \geq 5$.\\\\
 Now we can switch the role of $m$ and $n$ to complete the proof.
 \end{proof}
 \noindent
\begin{figure}[H]
  \centering
  \includegraphics[width=2in]{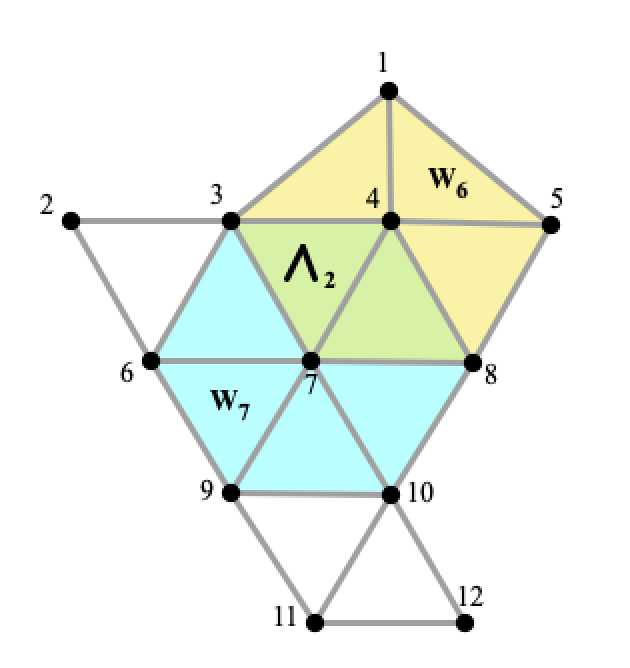}
  \captionof{figure}{Interlocking wheels as a subgraph of $G_F$}
  \label{fig:France_Wedge}
\end{figure}
\noindent
\textbf{Back to France}. Now we can compute $\chi(G_F,t)$.
First note that $W_6\wedge_2 W_7$ overlaps repeatedly with three 3-cycles in $K_2$ to produce $G_F$.
\begin{proof}[Chromatic Polynomial of $G_F$]
    \begin{align*}
      \chi(G_F,t)&=\frac{\chi(C_3,t)\chi(W_6 \wedge_2 W_7,t)\chi(C_3,t)\chi(C_3,t)}{\chi(K_2,t)\chi(K_2,t)\chi(K_2,t)}\\
      \\
&=\frac{\chi(C_3,t)^3\chi(\chi(W_6 \wedge_2 W_7,t)}{\chi(K_2,t)^3}\\ 
&=\frac{(t(t-1)(t-2))^3\frac{(t(t-2)(t-3)((t-2)^3-(-1)^5)((t-2)^4-(-1)^6)+t(t-2)^3((t-2)^2-(-1)^4)((t-2)^3-(-1)^5))}{(t-1)}}{(t(t-1))^3}\\
\\
&=(t-2)^3(t(t-2)(t-3)(t^2-5t+7)(t-1)(t-3)(t^2-4t+5)\\&+t(t-2)^3(t-1)(t-3)(t^2-5t+7)(t-1)\\
\\
&=t^{12}-23t^{11}+241t^{10}-1519t^9+6400t^8-18927t^7+40082t^6-60751t^5+64520t^4\\&-45656t^3+19328t^2-3696y .
    \end{align*}

\end{proof}
\noindent
Since $G_F$ contains $W_6$ as a subgraph, it is obvious it needs at least 4 colors to color it properly. This is confirmed by $\chi(G_F, t) = 0 $ for $t=1,2, \text{ and } 3$. However, we find that $\chi(G_F, 4) = 5184$. Thus the \\\\
\fbox{
\parbox{\textwidth}{
\textbf{number of ways to color the map of 12 contiguous regions of France is \\\\
\centerline{5184}.}
\\
}
}
3.\textbf{ America}.

\begin{figure}[H]
    \centering
    \includegraphics[width=4in]{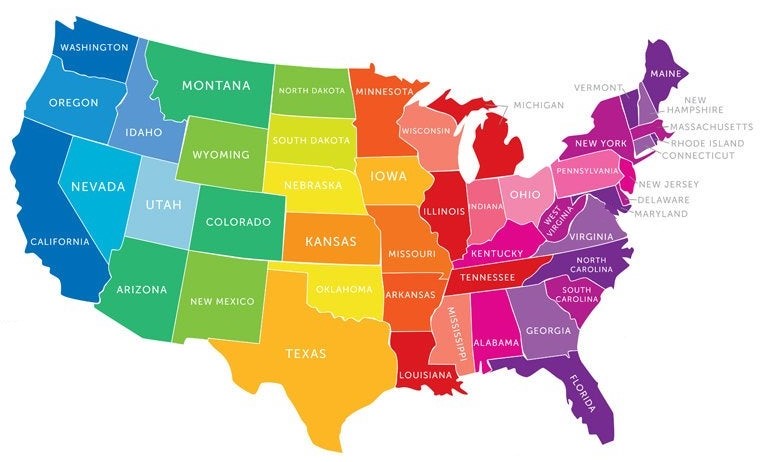}
    \caption{Map A, contains 48 states of the USA}
    \label{fig:map_USA}
\end{figure}

We associate the graph $G_A$ of the contiguous 48 States of America to the geographic map of the USA\footnote{Map of the US Retrived from: 
\url{https://img1.etsystatic.com/011/0/5380140/il_fullxfull.440914549_llhr.jpg}}, by placing a vertex on every state, putting an edge between two vertices if the corresponding states share a border (at more than just a corner). We do not include vertices for Hawaii and Alaska, because, it is enought to compute $\chi(G,t)$ for its connected components $G_j$ to get $\chi(G,t) = \prod_{j}\chi(G_j,t)$. \\

We attempted to calculate the chromatic polynomial for $G_A$. This was impossible to do entirely by hand, so we used modern computing technology to help with the calculations.

We tested different computing softwares, such as SAGE, Mathematica, MuPAD, and Maple. Each program was able to compute chromatic polynomials for relatively small graphs, but SAGE reached its limit when the number of vertices and edges was around 30 and 65 respectively. MuPAD couldn't handle a similarly sized graph either. Given a 39-vertex, 83-edge graph, Maple ran unsuccessfully for over an hour, while Mathematica was able to handle it in just 4 minutes. Mathematica, however, couldn't handle the entire $G_A$ with vertex-edge count of 48-105. From this, we knew our graph of the United States was too large for the best programs available to compute the chromatic polynomial. Because of its combinatorial nature, this is a problem of exponential complexity.

So, we combined the theory with our program of choice (Mathematica) to simplify our calculations.

Using the Chromatic Reduction Theorems, we were able to break down $G_A$ into subgraphs until it was reduced enough for Mathematica to compute the chromatic polynomials of these subgraphs. Using this approach, we calculated the chromatic polynomial of $G_A$ two different times in two different ways, each slightly differently.

\begin{figure}[H]
    \centering
    \includegraphics[width=3.2in]{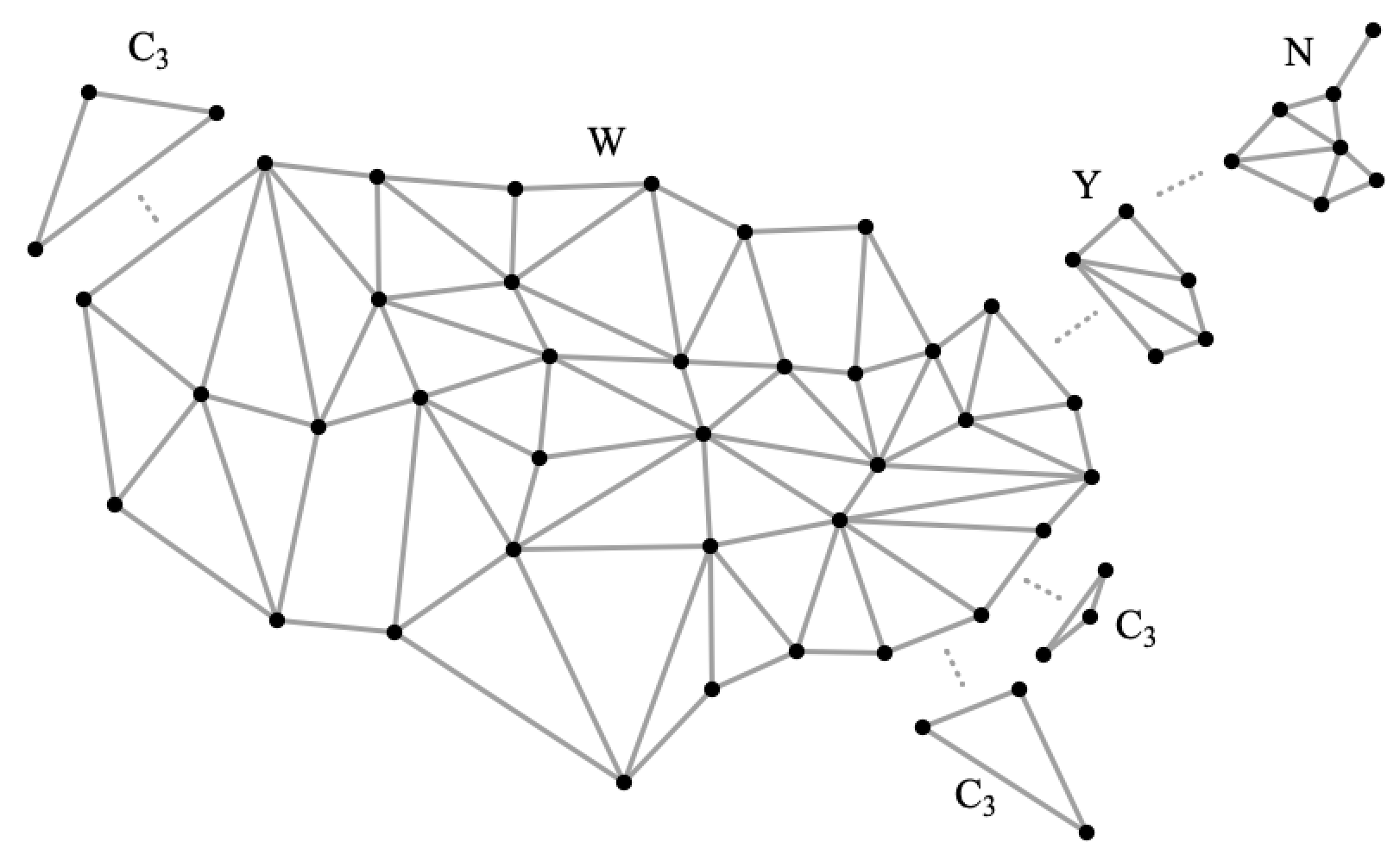}
    \caption{Attempt 1}
    \label{fig:attempt1}
\end{figure}

In the first attempt, we consider $G_A$ an overlap of $W$ and graphs $N$, $Y$, and 3 $C_3$'s(cf. Figure~\ref{fig:attempt1}). We let $N$ represent the northeastern-most states. We see that $Y$ overlaps $N$ and $W$, in $K_1$ and $K_2$ respectively. Each $C_3$ overlaps $W$ in $K_2$.

Using CRT-1,
\[\chi(Y,t)=\frac{\chi(C_3,t)^3}{\chi(K_2,t)^2}.\]

\noindent
We calculated $\chi(N,t)$ similarly.\\\\
\noindent
Finally, by a repeated application of CRT-1, we get\\\\ \begin{equation}\label{eq:4}
\chi(G_A,t)=\frac{\chi(W,t)\chi(Y,t)\chi(N,t)\chi(C_3,t)^3}{\chi(K_2,t)^4\chi(K_1,t)}
\end{equation}
\noindent
Then using Mathematica to compute $\chi(W, t)$, and plugging it in (\ref{eq:4}), we get\\\\
$\chi(G_A, t) =
 t^{48}-105t^{47}\newline+5404t^{46}-181689t^{45}\newline+4487296t^{44}-
86797239t^{43}\newline+1369003119t^{42}-18100363324t^{41}\newline+
204677484054t^{40}-2009741557171t^{39}\newline+17339117549604t^{38}-
132682763002081t^{37}\newline+907423360476887t^{36}-
5581169381630167t^{35}\newline+31031575427165032t^{34}-
156643500973559120t^{33}\newline+720455907112532420t^{32}-
3028205766124900090t^{31}\newline+11660587916045449786t^{30}-
41218581559720380212t^{29}\newline+133971262465065322950t^{28}-
400893632262902775367t^{27}\newline+1105490166090074304464t^{26}-
2811013847987117591939t^{25}\newline+6593222714427581969721t^{24}-
14265139633481475975539t^{23}\newline+28463416059570443463946t^{22}-
52346700555134790196556t^{21}\newline+88655379811509518107152t^{20}-
138105582619332483057236t^{19}\newline+197572598030181111248913t^{18}-
259057209010976705704331t^{17}\newline+310574028340459418761423t^{16}-
339434199963516594330980t^{15}\newline+336994442997740949240778t^{14}-
302626296141859413511498t^{13}\newline+244548011549333689537938t^{12}-
176716253459552869763068t^{11}\newline+113323961012010777670232t^{10}-
63883998375551455822512t^{9}\newline+31284178773965469140544t^{8}-
13106541029265128508800t^{7}\newline+4603738556785058047232t^{6}-
1318612298697138044928t^{5}\newline+295742443634384435712t^{4}-
48705983353199143936t^{3}\newline+5236695782665893888t^{2}-
275716645154500608t$\\\\
Evaluating it at $t=4$ gives $\chi(G_A, 4) = 12,811,729,152$. 
Thus,\\\\
\noindent
\fbox{
\parbox{\textwidth}{
\textbf{the number of ways to color the map of the contiguous 48 states of America is \\\\\centerline{12,811,729,152}}
\\
}
}
\section{Conclusion}
\begin{figure}[H]
    \centering
    \includegraphics[width=3.2in]{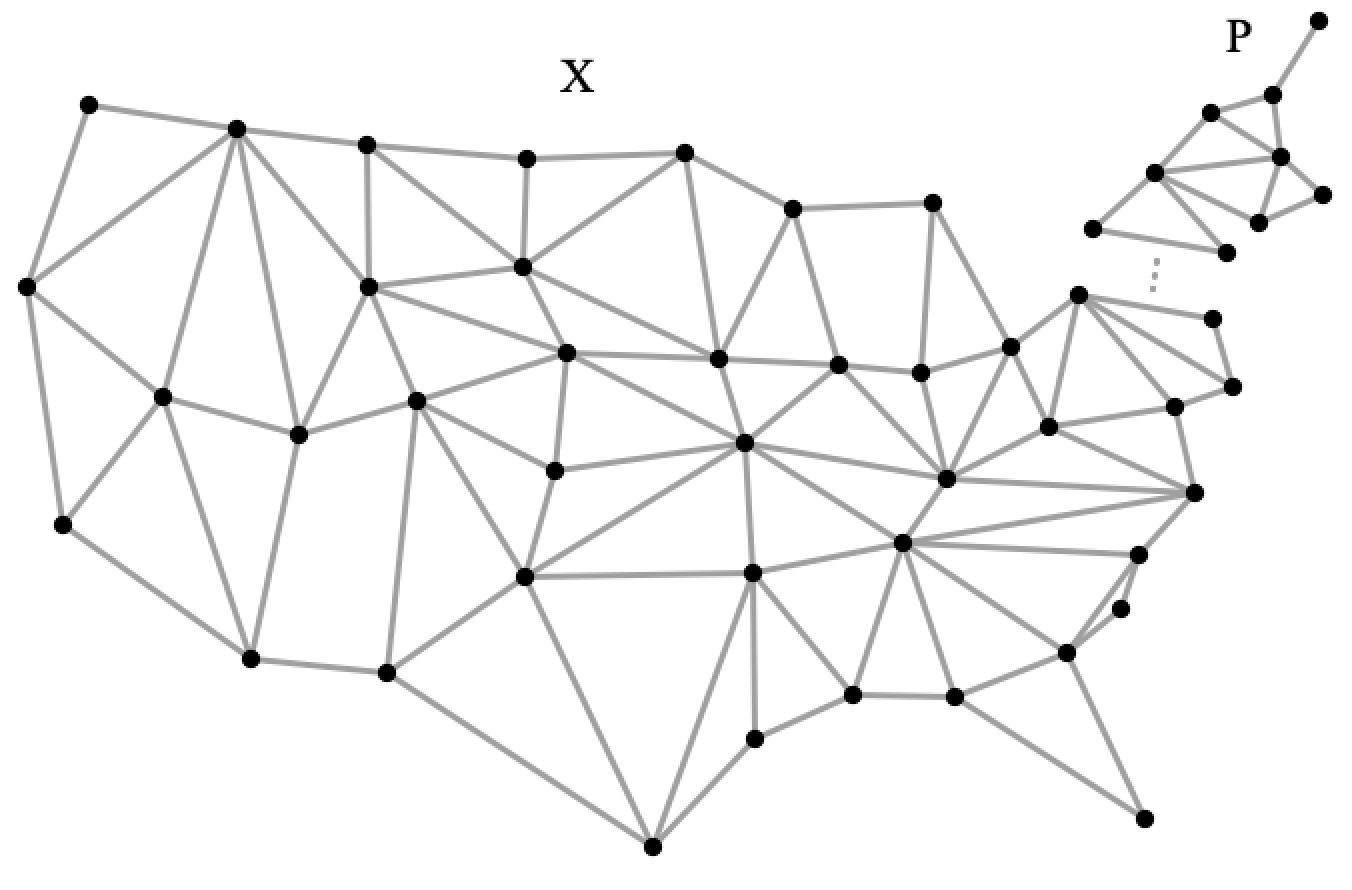}
    \caption{Attempt 2}
    \label{fig:attempt2}
\end{figure}

In Section~\ref{sec:coloringproblem}, we listed some properties of the chromatic polynomial of a graph. We now use them to check our answer for $\chi(G_A,t)$:
\begin{enumerate}[label=\arabic*)]
    \item It is true that $\chi(G_A,t)$ is monic, with integer coefficients.
        \item The degree of $\chi(G_A,t)$ is 48, the number of vertices in $G_A$.
    \item The coefficient of $t^{47}$ is $-105$, the number of edges of $G_A$ is 105.
    \item $G_A$ has no constant term, i.e. the constant term is zero.
    \item The degree of the lowest term in $\chi(G_A,t)$ is 1, which corresponds to the single component $G_A$ of the map of the contiguous 48 states of America.
    \item The coefficients of $\chi(G_A,t)$ clearly alternate between positive and negative.
        \item As Mathematica verified, the sum of the coefficients of $G_A$ is indeed zero.
\end{enumerate}
Thus, the polynomial we have found satisfies all of the above necessary conditions for a chromatic polynomial, and we are confident that our answer is correct.\\\\
For our second attempt, in order to double check the results of our first attempt, we consider $G_A$ an overlap of $X$ and $P$ in $K_2$, as pictured in Figure~\ref{fig:attempt2}.

\noindent
We see that \[\chi(G_A,t)=\frac{\chi(X,t)\chi(P,t)}{\chi(K_2,t)}\]
Surprisingly, Mathematica is able to compute the chromatic polynomial for X, a graph with 41 vertices and 93 edges. And of course, using methods previously discussed, $\chi(P,t)$ can be calculated by applying repeatedly the CRT - 1.

Again, it turns out that \[\chi(G_A,t)=\frac{\chi(X,t)\chi(P,t)}{\chi(K_2,t)}=\frac{\chi(W,t)\chi(Y,t)\chi(N,t)\chi(C_3,t)^3}{\chi(K_2,t)^4\chi(K_1,t)}\]
\noindent
As expected, $\chi(G,t) = 0$ for all $t\leq 3$. However, 
\[\chi(G_A,4)=12,811,591,729,152,\] the number of ways to color $G_A$.\\\\ 
\textbf{Remarks}.\\
1. To color the entire map of the USA, multiply it with $4 \cdot 4 = 16$ to include Alaska and Hawaii. If and when D.C becomes a state, the largest component $G_A$ of the map of the USA has to be enlarged and the new $\chi(G_A,t)$ computed again.

\section*{Acknowledgements}
The authors $\dagger$) would like to thank the College of Mathematical and Physical Sciences of BYU for financially supporting their research.\\\\

Rebekah Bassett, Department of Mathematics, BYU, Provo, UT 84602 USA email: rebekahbassett@yahoo.com \\

Jennifer Canizales, Department of Mathematics, BYU, Provo, UT 84602 USA email: jnydodge@gmail.com\\ 

Jasbar S. Chahal, Department of Mathematics, BYU, Provo, UT 84602 USA email: jasbir@math.byu.edu \\

Thomas Fackrell, Department of Mathematics, BYU, Provo, UT 84602 USA email: thomasfackrell@gmail.com \\

Vanessa Rico, Department of Mathematics, BYU, Provo , UT 84602 USA email: vrico@byu.edu

\end{document}